\documentclass[12pt]{amsart}
\renewcommand{\bar}{\overline}

\newcommand{\pa}{\partial}

\newcommand{\bb}{{\frac{\sqrt{-1}}{2\pi}}}

\newcommand{\ppp}[4]{\hspace{-3mm}
\begin{array}{l}
\begin{scriptscriptstyle}
{#1}\leq i\leq {#2}
\end{scriptscriptstyle}\\[-2mm]
\begin{scriptscriptstyle}
{#3}\leq j\leq {#4}
\end{scriptscriptstyle}
\end{array}
}

\newcommand{\pppp}[4]{\hspace{-2mm}
\begin{array}{l}
\begin{scriptscriptstyle}
{#1}\leq i\leq {#2}
\end{scriptscriptstyle}\\[-2mm]   
\begin{scriptscriptstyle}
{#3}\leq j\leq {#4}
\end{scriptscriptstyle}
\end{array}
}

\usepackage{amsmath,amsthm,amscd}
\newcommand{\frk}[1]{{\mathfrak{#1}}}

\title
[]{On the Futaki Invariants of 
Complete Intersections}
\author[]{Zhiqin Lu}
\date{June 3rd, 1998}
\address[Zhiqin Lu]
{Department of Mathematics\\
Columbia University\\
NY, NY 10027}
\email{lu@cpw.math.columbia.edu}

\newtheorem{theorem}{Theorem}[section]
\newtheorem{lemma}{Lemma}[section]
\newtheorem{cor}{Corollary}[section]
\newtheorem{prop}{Proposition}[section]

\newtheorem{ex}{Example}

\theoremstyle{remark}

\begin{document}
\maketitle

\numberwithin{equation}{section}

\newcommand{\kk}{K\"ahler~~}
\section{Introductions}
In 1983, Futaki~\cite{Fu1} introduced his invariants which 
generalize
the obstruction of Kazdan-Warner to prescribe Gauss curvature on $S^2$.
The Futaki invariants are defined for any compact K\"ahler manifold with
positive first Chern class that has nontrivial holomorphic vector fields.
Their vanishing are necessary conditions to the existence of
K\"ahler-Einstein metric on the underlying manifold.

Let $M$ be a compact K\"ahler manifold with positive
first Chern class $c_1(M)>0$. Choosing an arbitrary positive $(1,1)$ form 
$\omega$ in $c_1(M)$ as a K\"ahler metric on $M$, we can find a smooth
function $f$ on $M$, determined
up to a constant, such that the following
\begin{equation}\label{first}
Ric(\omega)-\omega=\bb\pa\bar\pa f
\end{equation}
holds.
Let ${\frk b}(M)$ be the Lie algebra of holomorphic vector fields on $M$.
The Futaki invariants are defined as
\[
F: {\frk b}(M)\rightarrow C,
F(X)=\int_M X(f)\omega^n.
\]

Ding and Tian~\cite{DT} took a further step in
introducing the Futaki invariants to Fano normal varieties.
This is not only a generalization of Futaki invariants to
singular varieties, but also has important application in
Kahler-Einstein geometry. In~\cite{T8}, the generalized K-energy
on normal varieties was introduced and a stability criteria
for the hypersurface or complete intersection was
established by using the notion of the generalized K-energy.
The Futaki invariants on singular varieties are related to
the stability of Fano manifolds due to the work of Tian
~\cite{T3}.  
To be more precise,
checking the $K$-stability of a Fano manifold is the same as checking 
the sign of the real part of the Futaki invariants on the 
 degenerations of the Fano manifold. Because of this,  we
need an effective way to compute the Futaki invariants
on singular varieties.

In this paper, we give a simple formula 
for the Futaki invariants
of Fano complete intersections.
The main theorem of this paper is:

\begin{theorem}\label{fud}
Let $M$ be the $N-s$ dimensional normal Fano variety in $CP^N$ defined
by the homogeneous polynomials $F_1,\cdots,F_s$ of degree
$d_1,\cdots,d_s$ respectively. Let $X$ be a holomorphic vector field on
$CP^N$ such that
\begin{equation}\label{01}
XF_i=\kappa_iF_i,\quad i=1,\cdots,s
\end{equation}
for constants $\kappa_1,\cdots,\kappa_s$.
Then the Futaki invariant $F(X)$ is
\begin{equation}\label{02}
F(X)=m^{N-s}\prod_{i=1}^s d_i
(-\sum^s_{i=1}\kappa_i+\frac{m}{N-s+1}
\sum_{i=1}^s\frac{\kappa_i}{d_i}),
\end{equation}
where $m=N+1-d_1-\cdots-d_s$.
\end{theorem}

\begin{cor}\label{fcor}
If $M$ is a hypersurface in $CP^N$ defined by the homogeneous polynomial
$F$ of degree $d$ and if $XF=\kappa F$, then
\[
F(X)=-(N+1-d)^{N-1}\frac{(N+1)(d-1)}{N}\kappa.
\]
In particular, $Re\, F(X)$ and $ -Re\,\kappa$ have the same sign.
\end{cor}

The formula is new even in the case when $M$ is a hypersurface or an
orbifold.
If the zero locus of the holomorphic vector field
$X$ on $M$ is a smooth manifold,
then
using the residue formula of the Atiyah-Bott-Lefschetz type, Futaki was
able
to
 develop a method to compute his invariants by the information of the
vector field $X$ and the manifold $M$ near the zero locus of the vector
field~\cite{Fu2}.
In~\cite{DT}, the authors developed the  method to compute the Futaki
invariants on 2-dimensional K\"ahler orbifolds.  In~\cite{Ma1}, the Futaki
invariants
for toric varieties were calculated.

{\bf Acknowledgment.} 
The author thanks G. Tian for his mathematical insights and
encouragements  during the preparation of this
paper. He also thanks D. Phong, M.
Kuranishi, H. Pinkham, L. Borisov and Z. Wu for the discussions of this
topic.

\section{Preliminaries}

Let $X=(a_{ij})_{\pppp{1}{N}{1}{N}}\in sl(N+1,C)$ be
a traceless $(N+1)\times (N+1)$ matrix. $X$ 
defines a holomorphic vector field
\begin{equation}\label{XX}
X=\sum_{i,j=0}^N a_{ij}z^i\frac{\pa}{\pa z^j}
\end{equation}
on $C^{N+1}$ and a smooth function 
\begin{equation}\label{cc}
\theta=-X\log\sum_{i=0}^N |z^i|^2
\end{equation}
on $C^{N+1}-\{0\}$, 
where $(z^0,\cdots,z^N)$ is the coordinates of $C^{N+1}$.

Both $X$ and $\theta$ descends to a vector field and 
a smooth function on
the 
projective space $CP^N$, respectively.
Let $\omega_{FS}$ be the Fubini-Study
metric of $CP^N$. 
Then we have the following relation:
\begin{equation}\label{21}
i(X)\omega_{FS}=-\bb\bar\pa\theta.
\end{equation}

Suppose  $M$ is an $n$-dimensional $Q$-Fano normal variety in $CP^N$ and
suppose that
$\omega=\alpha\omega_{FS}\in c_1(M)$ for a constant $\alpha>0$. If $X$ is
a vector field on $CP^N$ such that the one parameter group generated 
by the real part of $X$ leaves $M$ invariant,
we say that $X$ is tangent to $M$.
Suppose $(x^1,\cdots,x^n)$ is the local holomorphic coordinates at some
smooth point $p$ of $M$.  
Equation~\eqref{21} can be written as
\begin{equation}\label{22}
X^i=-\alpha g^{i\bar j}\frac{\pa \theta}{\pa \bar x^i},\quad
i=1,\cdots,n,
\quad X=\sum_{i=1}^n X^i\frac{\pa}{\pa x^i},
\end{equation}
where $(g^{i\bar j})$ is the inverse matrix of $(g_{i\bar j})$
and $(g_{i\bar j})$ is the metric matrix of $\omega$.

We
define the divergence of $X$ on $M$ by
\begin{equation}\label{x233}
div\, X=\frac{\pa X^i}{\pa x^i} +\sum_{i=1}^n X^i
\frac{\pa}{\pa x^i}\log\det (g_{i\bar j}).
\end{equation}

The following lemma is the observation on which the whole paper is based.

\begin{lemma}[\cite{T3}]\label{TT3}
 If $M$ is a normal projective variety, then
\begin{equation}\label{xx}
div\, X-\alpha\theta+X(f)=const
\end{equation}
where the function $f$ is defined as
\begin{equation}\label{23}
Ric(\omega)-\omega=\bb\pa\bar\pa f
\end{equation}
\end{lemma}

{\bf Proof:} A straightforward computation yields
\[
\bb\bar\pa\, div\, X=-i(X) Ric(\omega).
\]
By Equation~\eqref{23}, we see that
\[
\bb\bar\pa\, div\, X=-i(X)\omega-\bb i(X)\pa\bar\pa f
=\bb \alpha\bar\pa\theta-\bb\bar\pa X(f).
\]
Thus $div\, X-\alpha\theta+X(f)$ is a holomorphic function on 
the normal variety $M$ which  must be a constant.

\qed

\begin{cor}\label{cor1}
The Futaki invariant can be written as
\[
F(X)= const\,vol(M)+\alpha\int_M\theta\omega^n.
\]
\end{cor}

\qed

\section{An Explicit Expression of the Function $f$}
Suppose $M$ 
 is a complete intersection
of $CP^N$. That is, $M$ is the zero locus of 
homogeneous polynomials
$F_1,\cdots, F_s$ in $CP^N$
with degree $d_1,\cdots,d_s$,
respectively and the dimension of $M$ is  $N-s$. 
By the adjunction formula, the
anticanonical bundle of $M$ is
\[
K_M^{-1} = (N+1-d_1-\cdots -d_s)H,
\]
where $H$ is the hyperplane bundle of $CP^N$. 
We assume that $M$ is a normal variety.
There is a unique function $f$ (up to a
constant), defined  on the regular part $M_{reg}$ of $M$,
such that if $\omega=(N+1-d_1-\cdots -d_s)\omega_{FS}|_M$, then
\[
Ric (\omega)-\omega=\bb\pa\bar\pa f\qquad\text{on\,}M_{reg}.
\]

In this section, we write out the above function $f$ explicitly.
The idea is to trace the proof of the well known adjunction formula.
But here we work on the metric level rather than the cohomological
level. This makes the notations a little bit complicated.

\label{se}
We begin by the
following general setting: Let $V$ be a K\"ahler manifold of dimension $N$
and $V_1$ be a submanifold of dimensional  $N-s$ defined  by
holomorphic
functions $G_1=G_2=\cdots =G_s=0$. Suppose $U_1$ is an open set of $V$
such
that
\begin{enumerate}
\item
$(x^1,\cdots,x^N)$ is a local holomorphic coordinate system of $V$;
\item
On $U_1$, we have 
\[
rank\frac{\pa (G_1,\cdots,G_s)}{\pa (x^1,\cdots,x^s)}=s;
\]
\item 
There are holomorphic functions $f_1,\cdots,f_s$ on $U_1\cap V_1$ such
that
\[
\left\{
\begin{array}{l}
x^1=f_1(x^{s+1},\cdots,x^N)\\
x^2=f_2(x^{s+1},\cdots,x^N)\\
\cdots\cdots\\
x^s=f_s(x^{s+1},\cdots,x^N)
\end{array}
\right..
\]
\end{enumerate}
In particular, $(x^{s+1},\cdots,x^N)$ is the local holomorphic coordinate
system
of
$U_1\cap V_1$.

Suppose $g_1=\sum_{\pppp{s+1}{N}{s+1}{N}} g_{1i\bar j} dx^i\wedge
d\bar{x}^j$ is the restriction of the K\"ahler metric of $V$
on $U_1\cap V_1$. Define
$\det g_1=\det (g_{1i\bar j})_{\ppp{s+1}{N}{s+1}{N}}$.
Of course, $\det g_1$ is not a global function
on $V_1$. In order to 
study the change of the $\det g_1$ with respect to the change of the
local holomorphic coordinates,
we assume that
there is another neighborhood 
 $(U_2, (y^1,\cdots,y^N))$ 
of $V$ such that $U_1\cap U_2\cap V_1\neq\emptyset$. As before, we assume
that
\[
rank \frac{\pa (G_1,\cdots, G_s)}{\pa (y^1,\cdots,y^s)}=s
\]
and on $U_2\cap V_1$, we have
\[
\left\{
\begin{array}{l}
y^1=g_1(y^{s+1},\cdots,y^N)\\
y^2=g_2(y^{s+1},\cdots,y^N)\\
\cdots\cdots\\
y^s=g_s(y^{s+1},\cdots,y^N)\\
\end{array}
\right.
\]
for holomorphic functions $g_1,\cdots,g_s$ on $U_2\cap V_1$.
$(y^{s+1},\cdots,y^N)$ is the local holomorphic coordinate
system  of
$U_2\cap V_1$. Let 
\[
g_2=\sum_{\pppp{s+1}{N}{s+1}{N}}g_{2i\bar j}
dy^i\wedge
d\bar y^j
\]
be the restriction of the K\"ahler metric of $V$
on $U_2\cap V_1$. Define
$\det g_2=\det (g_{2i\bar j})_{\ppp{s+1}{N}{s+1}{N}}$.
Then we have

\begin{prop}[Adjunction Formula]\label{p41}
With the above notations, \linebreak on $V_1\cap U_1\cap
U_2\neq\emptyset$, we have
\[
\det (g_1) =\det (g_2)
\left |\det (\frac{\pa y^j}{\pa x^i})_{\ppp{1}{N}{1}{N}}\right |^2\cdot
\frac{\displaystyle
\left |\det\frac{\pa (G_1,\cdots, G_s)}{\pa (y^1,\cdots, y^s)}\right
|^2}
{\displaystyle \left
|\det\frac{\pa (G_1,\cdots, G_s)}{\pa (x^1,\cdots, x^s)}\right |^2}.
\]
\end{prop}
{\bf Proof:}
Let
\[
\left\{
\begin{array}{ll}
\tilde x^i=x^i-f_i(x^{s+1},\cdots,x^N) & 1\leq i\leq s\\
\tilde x^i=x^i & i>s
\end{array}
\right.,
\]
and
\[
\left\{
\begin{array}{ll}
\tilde{y}^j=y^j-g_j(y^{s+1},\cdots,y^N) & 1\leq j\leq s\\
\tilde{y}^j=y^j & j>s
\end{array}
\right..
\]
Then $V_1$ is locally defined by $\tilde x^i=0\, (1\leq i\leq s)$
or $\tilde y^j=0\, (1\leq j\leq s)$. In particular, on $V_1$ we have
\begin{equation}\label{pss}
\frac{\pa\tilde x^i}{\pa\tilde y^j}=0,\quad 1\leq i \leq  s,\quad s+1\leq
j\leq N.
\end{equation}

Before going further, we make the following conventions:
\begin{enumerate}
\item 
$\det (\frac{\pa x^i}{\pa y^j})_{\pppp{1}{N}{1}{N}}$
and 
$\det (\frac{\pa\tilde x^i}{\pa\tilde y^j})_{\pppp{1}{N}{1}{N}}$
denote the restriction to $V_1$ of the Jocobi determinant
of the transform
$(y^1,\cdots,y^N)\rightarrow (x^1,\cdots,x^N)$
and 
$(\tilde y^1,\cdots,\tilde y^N)\rightarrow (\tilde x^1,\cdots,\tilde x^N)$
respectively;
\item
$\det (\frac{\pa\tilde x^i}{\pa\tilde y^j})_{\pppp{1}{s}{1}{s}}$
and
$\det (\frac{\pa\tilde x^i}{\pa\tilde y^j})_{\pppp{s+1}{N}{s+1}{N}}$
are the determinant of the submatrices of
$(\frac{\pa\tilde x^i}{\pa\tilde y^j})_{\pppp{1}{N}{1}{N}}$;
\item
Since $(x^{s+1},\cdots,x^N)$ and $(y^{s+1},\cdots,y^N)$
 are local coordinates of $U_1\cap V_1$ and $U_2\cap V_2$
respectively, there is the transform \linebreak
$(y^{s+1},\cdots,y^N)\rightarrow (x^{s+1},\cdots,x^N)$
by
\[
x^i=x^i(g_1(y^{s+1},\cdots,y^N),\cdots,g_s(y^{s+1},\cdots,y^N),
y^{s+1},\cdots,y^N)
\]
for $\quad s+1\leq i\leq N$
and $\det(\frac{\pa x^i}{\pa y^j})_{\pppp{s+1}{N}{s+1}{N}}$
denotes the Jacobi determinant of the above
transform.
\end{enumerate}

If $s+1\leq i,j\leq N$, then
\[
\frac{\pa\tilde x^i}{\pa\tilde y^j}
=
\frac{\pa x^i(g_1,\cdots,g_s,y^{s+1},\cdots,y^N)}{\pa y^j}.
\]

Using  Equation~\eqref{pss},  we have
\begin{equation}\label{422}
\det (\frac{\pa x^i}{\pa y^j})_{\ppp{1}{N}{1}{N}}
=
\det (\frac{\pa \tilde x^i}{\pa \tilde y^j})_{\ppp{1}{s}{1}{s}}
\cdot
\det (\frac{\pa x^i}{\pa y^j})_{\ppp{s+1}{N}{s+1}{N}}.
\end{equation}

If $1\leq i, j\leq s$, then
\[
\frac{\pa\tilde x^i}{\pa\tilde y^j}
=
\sum_{k=1}^{s} \frac{\pa \tilde x^i}{\pa y^k}\cdot
\frac{\pa y^k}{\pa\tilde y^j}
=
\frac{\pa x^i}{\pa y^j}
-
\sum_{k=s+1}^N\frac{\pa f_i}{\pa x^k}\cdot
\frac{\pa x^k(g_1,\cdots,g_s,y^{s+1},\cdots,y^N)}{\pa y^j}.
\]
Thus we have
\begin{equation}\label{pl}
\sum_{i=1}^{s}\frac{\pa G_k}{\pa x^i}\cdot
\frac{\pa\tilde x^i}{\pa\tilde y^j}
=
\frac{\pa G_k}{\pa y^j}.
\end{equation}
The proposition follows from Equation~\eqref{422} and ~\eqref{pl}.

\qed

We are going to use the above proposition in the case of complete
intersections of $CP^N$.
Since $M$ is defined by the zero locus of homogeneous functions, we must
make some necessary adjustment because homogeneous polynomials are
\textsl{not} functions on $CP^N$.

Let $U_\alpha (\alpha=0,\cdots,n)$ be the standard covering
of $CP^N$ defined
by $U_\alpha=\{ Z_\alpha\neq 0\}$ where $[Z_0,\cdots,Z_N]$ is the 
homogeneous coordinates of $CP^N$. Suppose
$z_i^\alpha={Z_i}/{Z_\alpha}, (i\neq\alpha)$
be the standard coordinates on $U_\alpha$. Let
\[
F_i^\alpha(z_0^\alpha,\cdots,\hat{z_\alpha^\alpha},\cdots,z_N^\alpha)
=F(Z_0,\cdots,Z_N)/(Z_\alpha)^{d_i},\quad 1\leq i\leq s.
\]
For each $\{\alpha,\beta_1,\cdots,\beta_s\}\subset\{0,\cdots,N\}$, define
\[
U_{\alpha\beta_1\cdots\beta_s}
=
\{
p\in U_\alpha |
rank \frac{\pa (F_1^\alpha,\cdots,F_s^\alpha)}
{\pa (z^\alpha_{\beta_1},\cdots,z^\alpha_{\beta_s})}
=s\}.
\]
Then it is clear that
$\cup U_{\alpha\beta_1\cdots\beta_s}\supset M_{reg}$.

At each point $p\in M_{reg}\cap U_{\alpha\beta_1\cdots\beta_s}$, 
$(z^\alpha_\beta, \beta\notin\{\alpha,\beta_1,\cdots,\beta_s\})$ can be
used as local coordinate system at $p$. Let
$g_{\alpha\beta_1\cdots\beta_s}$ be the corresponding metric matrix and let 
$\det g_{\alpha\beta_1\cdots\beta_s}$ be its determinant.
Define
\begin{equation}\label{4aa}
\xi_{\alpha\beta_1\cdots\beta_s}
=
\det (g_{\alpha\beta_1\cdots\beta_s})
\left |\det\frac{\pa (F_1^\alpha,\cdots,F_s^\alpha)}
{\pa (z^\alpha_{\beta_1},\cdots, z_{\beta_s}^\alpha)}\right |^2
(1+\sum_{i\neq \alpha} |z_i^\alpha|^2)^{N+1-d_1-\cdots-d_s}.
\end{equation}
Then we have
\begin{lemma}
$\{\xi_{\alpha\beta_1\cdots\beta_s}\}$ defines a global positive function
of $M_{reg}$.
\end{lemma}

{\bf Proof:}
A straightforward computation shows (cf.~\cite[pp. 146]{GH})
\[
dz^{\alpha_2}_0\wedge\cdots \widehat{dz^{\alpha_2}_{\alpha_2}}
\cdots\wedge dz^{\alpha_2}_N
=\left(\frac{Z_{\alpha_1}}{Z_{\alpha_2}}\right)^{N+1}
dz^{\alpha_1}_0\wedge\cdots \widehat{dz^{\alpha_1}_{\alpha_1}}
\cdots\wedge dz^{\alpha_1}_N.
\]
The lemma follows from 
Proposition~\ref{p41},
Equation~\eqref{4aa} and the above equation.

\qed

\begin{theorem}\label{thm411}
Let $f$ be the function on $M_{reg}$, defined by
\begin{equation}\label{ps}
f=-\log \xi_{\alpha\beta_1\cdots\beta_s}\quad on\quad
U_{\alpha\beta_1\cdots\beta_s}.
\end{equation}
Then
\[
\bb\pa\bar\pa f= Ric(\omega)-\omega,
\]
where $\omega=(N+1-d_1-\cdots-d_s)\omega_{FS}|_M$.
\end{theorem}

\qed

\newcommand{\vv}{{\alpha\beta_1\cdots\beta_s}}
\newcommand{\vvv}{{\alpha,\beta_1,\cdots,	\beta_s}}

\newcommand{\ttt}{\frac{\pa (F_1^\alpha,\cdots,F^\alpha_s)}
{\pa (z^\alpha_{\beta_1},\cdots,z^\alpha_{\beta_s})}}
\section{The Trace of the Action on $F_1,\cdots,F_s$}
Let $M$ be the  variety defined in the previous section.
The vector field $X=\sum_{i,j=0}^N a_{ij} Z_j\frac{\pa}{\pa Z_i}$
naturally acts on $F_i$ by
\begin{equation}\label{5pq}
XF_i=\sum_{i,j=0}^N a_{ij} Z_j\frac{\pa F_i}{\pa Z_i},\quad
i=1,\cdots,s.
\end{equation}
Suppose ${\mathcal V}$ is the vector space spanned by $F_1,\cdots,F_s$.
Since $X$ is tangent to $M$, $X$ is an automorphism on ${\mathcal V}$.

The main result of this section is,
\begin{theorem}\label{thmu}
Let $\kappa$ be the trace of the automorphism of $X$ on ${\mathcal V}$.
Then
\[
div\,X+X(f)-(N+1-d_1\cdots-d_s)\theta=-\kappa.
\]
\end{theorem}

{\bf Proof:} We adopt all  notations from last section.
Consider a smooth point $p$ of $M_{reg}$ in
$U_{\alpha\beta_1\cdots\beta_s}$.
From Equation~\eqref{cc}, the function $\theta$ 
in the local coordinates $\{z_j^\alpha, j\neq\alpha\}$ is
\begin{equation}\label{abd}
\theta=-\tilde X\log (1+\sum_{i\neq\alpha}|z_i^\alpha|^2)-a_{\alpha\alpha}
-\sum_{j\neq\alpha} a_{\alpha j}z_j^\alpha.
\end{equation}
By~\eqref{XX}, let
\begin{equation}\label{aedd}
\tilde X=\sum_{i\neq\alpha}((\sum_{j\neq\alpha}
a_{ij}z_j^\alpha-a_{\alpha j}z_j^\alpha z_i^\alpha)
+(a_{i\alpha}-a_{\alpha\alpha}z_i^\alpha))
\frac{\pa}{\pa z_i^\alpha}.
\end{equation}
Let the holomorphic vector field $X$ on $M$ be written as
\begin{equation}\label{abe}
 X=\sum_{i\notin\{\vvv\}} X^i\frac{\pa}{\pa z_i^\alpha}.
\end{equation}
If $i: M\rightarrow CP^N$ is the embedding, then $i_*X=\tilde X$.

By Equation~\eqref{4aa}, ~\eqref{ps}, ~\eqref{abd}, ~\eqref{abe} 
and the definition of $div\, X$ in Equation ~\eqref{x233}, we see that
\begin{align}\label{5a}
\begin{split}
&div\,X+X(f)-(N+1-d_1-\cdots-d_s)\theta\\
&=\sum_{i\notin\{\vvv\}}\frac{\pa X^i}{\pa z_i^\alpha}
-X\log\det\ttt\\
&\qquad +(N+1-d_1-\cdots-d_s) 
(a_{\alpha\alpha}+\sum_{j\neq\alpha} a_{\alpha j}z_j^\alpha).
\end{split}
\end{align}

Before going on, we need  a general elementary lemma. To begin, we use
the general setting on page~\pageref{se}. In addition, we let $X$ be a
holomorphic vector field of $V$ such that $X$ is tangent to $V_1$.
In what follows, we  temporary distinguish the $X$ on $V$ and the
$X$ on $V_1$. So let's denote the $X$ on $V$ to be $\tilde X$. In the
local coordinates, $\tilde X$ is
\[
\tilde X=\sum_{i=1}^N \tilde X^i\frac{\pa}{\pa x^i}.
\]
Then $X$ on $V_1$ can be written as
\[
 X=\sum^N_{i=s+1} X^i\frac{\pa}{\pa x^i},
\]
where from the chain rule,
\begin{equation}\label{51}
 X^i=\tilde X^i(f_1,\cdots,f_s,X^{s+1},\cdots,X^N),\quad
s+1\leq i\leq N.
\end{equation}
If $i:V_1\rightarrow V$ is the embedding, then
$i_* X=\tilde X$.

We have the following elementary lemma:

\begin{lemma}\label{opq}
Let
\[
A=\frac{\pa(G_1,\cdots,G_s)}{
\pa(x^1,\cdots,x^s)}
\]
and let $A_{ij}=\frac{\pa G_i}{\pa x^j}$
for $i,j=1,\cdots s$ and $(A^{ij})$ be the inverse matrix
of $A$. Then on $V_1$, we have
\[
\sum_{i=s+1}^N\frac{\pa X^i}{\pa x^i}
-\tilde X\log\det A
=
\sum_{i=1}^N\frac{\pa \tilde X^i}{\pa x^i}
-\sum_{i,j=1}^s A^{ji}\frac{\pa}{\pa x^j} \tilde XG_i.
\]
\end{lemma}

{\bf Proof:} By definition,
\[
\tilde X\log\det A=\sum_{k=1}^N\sum_{i,j=1}^s
A^{ji} \tilde X^k\frac{\pa^2 G_i}{\pa x^j\pa x^k}.
\]
We can write the above equation as
\begin{equation}\label{477}
\tilde X\log\det A=\sum_{i,j=1}^s(A^{ji}\frac{\pa}{\pa x^j} 
\tilde  XG_i
-\sum_{k=1}^N A^{ji}
\frac{\pa\tilde  X^k}{\pa x^j}
\frac{\pa G_i}{\pa x^k}).
\end{equation}
By the implicit differentiation, we see that on $V_1$,
\begin{equation}\label{488}
\frac{\pa f_j}{\pa x^k}
=-\sum^s_{i=1} A^{ji}
\frac{\pa G_i}{\pa x^k},
\quad,
j=1,\cdots s, k=s+1,\cdots,N.
\end{equation}
Using Equation~\eqref{477} and~\eqref{488},
\[
\tilde X\log\det A=\sum_{i,j=1}^s A^{ji}\frac{\pa}{\pa x^j} \tilde XG_i
+\sum_{k=s+1}^N\sum_{j=1}^s
\frac{\pa\tilde X^k}{\pa x^j}
\frac{\pa f_j}{\pa x^k}
-\sum_{i=1}^s\frac{\pa\tilde X^k}{\pa x^i}.
\]
The lemma follows from the above identity and the fact that
\[
\sum_{i=s+1}^N\frac{\pa X^i}{\pa x^i}
=\sum_{i=s+1}^N\frac{\pa\tilde X^i}{\pa x^i}
+\sum_{i=s+1}^N\sum_{j=1}^s\frac{\pa\tilde X^i}{\pa x^j}\cdot
\frac{\pa f_j}{\pa x^i}.
\]

\qed

Go back to the proof of the theorem. Let
\[
A=\ttt
\]
and temporary denote $\tilde X$ to be the vector field $X$ on $U_\alpha
=\{Z_\alpha\neq 0\}$.
The representation of $\tilde X$ is in Equation~\eqref{aedd}.
Obviously
\[
X\log\det\ttt=\tilde X\log\det \ttt.
\]
Using the Lemma~\ref{opq}, Equation~\eqref{5a} becomes
\begin{align}\label{4d7}
\begin{split}
&div\, X+X(f)-(N+1-d_1-\cdots-d_s)\theta\\
&=\sum_{i\neq\alpha}\frac{\pa\tilde  X^i}{\pa z_i^\alpha}
-\sum_{i=1}^s\sum_{j\in\{\beta_1,\cdots,\beta_s\}}A^{ji}\frac{\pa}{\pa
z_j^\alpha} \tilde XF_i^\alpha\\
&\qquad +(N+1-d_1-\cdots-d_s)(a_{\alpha\alpha}
+\sum_{j\neq\alpha} a_{\alpha j} z_j^\alpha).
\end{split}
\end{align}

Since $\sum a_{ii}=0$, a simple calculation gives
\begin{equation}\label{m22}
\sum_{i\neq\alpha}\frac{\pa\tilde  X^i}{\pa z_i^\alpha}
=
-(N+1)(a_{\alpha\alpha}+\sum_{j\neq\alpha} a_{\alpha j}z_j^\alpha).
\end{equation}
Recall the definition of  $XF_i$ in Equation~\eqref{5pq}. We
see that
for $i=1,\cdots,s$,
\[
\tilde XF_i^\alpha=X\frac{F_i}{Z_\alpha^{d_i}}
=\frac{XF_i}{Z_\alpha^{d_i}}
-d_i\frac{F_i}{Z_\alpha^{d_i}}
(a_{\alpha\alpha}+\sum_{j\neq\alpha} a_{\alpha j}z_j^\alpha).
\]
Thus on $M_{reg}$
\begin{align}\label{m23}
\begin{split}
&\qquad\qquad\sum_{i=1}^s\sum_{j\in\{\beta_1,\cdots,\beta_s\}}
A^{ji}\frac{\pa}{\pa z_j^\alpha}\tilde X F_i^\alpha\\
&=\sum_{i=1}^s\sum_{j\in\{\beta_1,\cdots,\beta_s\}}
A^{ji}Z_\alpha
\left(\frac{{\frac{\pa}{\pa Z_j} XF_i}}
{{Z_\alpha^{d_i}}}
-d_i\frac{{\frac{\pa F_i}{\pa Z_j}}}
{Z_\alpha^{d_i}}(a_{\alpha\alpha}
+\sum_{j\neq\alpha} a_{\alpha j}z_j^\alpha)\right)\\
&=\kappa-(\sum_{i=1}^s d_i)
(a_{\alpha\alpha}+\sum_{j\neq\alpha}z_j^\alpha),
\end{split}
\end{align}
where we used the fact that on $M_{reg}$, $F_1=\cdots=F_s=0$ in the 
second identity.
The theorem follows from Equation~\eqref{4d7}, 
\eqref{m22} and \eqref{m23}.

\section{The Computation of the $\theta$ Invariants}
Let $M$ be the complete intersection defined in \S 3.
Let
$M_0=CP^N$. $M_k=N_1\cap\cdots\cap N_k, (k=1,\cdots,s)$. 
Then $M_s=M$.

We
 assume that
$XF_i=\kappa_iF_i,\quad i=1,\cdots,s$.

Let $[Z_0,\cdots,Z_N]$ be the homogeneous coordinates
of $CP^N$. Define
\[
\xi_i=\frac{|F_i|^2}{\sum_{i=0}^N(|Z_i|^2)^{d_i}},
\quad, i=1,\cdots,s.
\]
Then $\xi_i$'s are global smooth functions on $CP^N$. 

In this section, we compute the $\theta$-invariant
$\int_M\theta\omega^n$, where 
for simplifying the notations, we assume  that
$\omega=\omega_{FS}$ is the Fubini-Study metric of the $CP^N$. The key
result is the
following:

\begin{lemma}\label{p51}
For $k=2,\cdots,s$, we have
\begin{equation}\label{6a}
\int_{M_k}(\theta+\omega)^{N-k+1}
=d_k\int_{M_{k-1}}(\theta+\omega)^{N-k+2}
+\kappa_k d_1\cdots d_{k-1},
\end{equation}
and in addition, we have
\begin{equation}\label{6b}
\int_{M_1}(\theta+\omega)^N=\kappa_1.
\end{equation}
\end{lemma}

{\bf Proof:} We have the following identities for $k=1,\cdots,s$
\begin{align}\label{6c}
\begin{split}
&\quad\bar\pa(\bb\pa\log\xi_k\wedge\theta\omega^{N-k})
+i(X)(\frac{\pa\log\xi_k}{N-k+1}\omega^{N-k+1})\\
&=-\bb\pa\bar\pa\log\xi_k\wedge\theta\omega^{N-k}
+\frac{1}{N-k+1}X\log\xi_k\wedge\omega^{N-k+1}.
\end{split}
\end{align}
Integration against $M_{k-1}$ gives
\[
\int_{M_{k-1}}\bb\pa\bar\pa\log\xi_k\wedge\theta\omega^{N-k}
=\frac{1}{N-k+1}
\int_{M_{k-1}} X\log\xi_k\wedge\omega^{N-k+1}.
\]

Since for $k=1,\cdots,s$,
\[
\begin{array}{c}
\bb\pa\bar\pa\log\xi_k=[N_k]-d_k\omega\\
X\log\xi_k=\kappa_k+d_k\theta,
\end{array}
\]
where $[N_k]$ is the divisor of the zero locus of $F_k$,
we have
\[
\int_{M_k}\theta\omega^{N-k}-d_k\int_{M_{k-1}}\theta\omega^{N-k+1}
=\frac{1}{N-k+1}\int_{M_{k-1}}(\kappa_k+d_k\theta)\omega^{N-k+1}.
\]
Thus
\begin{align*}
&\int_{M_k}\theta\omega^{N-k}
=\frac{N-k+2}{N-k+1} d_k\int_{M_{k-1}}\theta\omega^{N-k+1}
+\frac{\kappa_k d_1\cdots d_{k-1}}{N-k+1}.
\end{align*}

So ~\eqref{6a} is 
proved. To prove ~\eqref{6b}, let's first see that by ~\eqref{6c},
\[
\int_{M_1}\theta\omega^{N-1}=\frac{N+1}{N}d_1
\int_{CP^N}\theta\omega^N+\frac{\kappa_1}{N}.
\]

Then ~\eqref{6b} follows from the simple fact that
\[
\int_{CP^N}\theta\omega^N=0.
\]

\qed

Equation~\eqref{6a} can be rewritten as 
\[
\frac{1}{d_1\cdots d_k}\int_{M_k} (\theta+\omega)^{N-k+1}
=\frac{1}{d_1\cdots d_{k-1}}\int_{M_{k-1}}
(\theta+\omega)^{N-k+2}
+\frac{\kappa_k}{d_k}
\]
for $k=2,\cdots,s$. Thus we have

\begin{theorem}\label{py}
 With the notations as above, we have
\[
\int_M\theta\omega^{N-s}=
\frac{d_1\cdots d_s}{N-s+1}
\sum_{k=1}^s\frac{\kappa_k}{d_k}.
\]
\end{theorem}

\qed

Now we  prove the main theorem of this paper:

{\bf Proof of Theorem~\ref{fud}:} The theorem follows from
Theorem~\ref{thmu}, Theorem
~\ref{py} and the fact that
\[
\omega=(N+1-d_1-\cdots-d_s)\omega_{FS}\in c_1(M).
\]

\qed

\section{Examples}
In this section, we 
use our formula to 
compute some examples given by Ding-Tian~\cite{DT}, Jeffres~\cite{JE}
and Wu~\cite{Wu}. 
Recall that, the Futaki invariant defined in ~\cite{DT} and ~\cite{JE}
differ from us by a factor 3 in the case of surfaces.
So in what follows, the Futaki invariant $F(X)$ is actually three times
the Futaki invariant in the previous sections.

\begin{cor}\label{cor61}
With the notations as in Theorem~\ref{fud},
if $M$ is the cubic surface in $CP^3$, then the Futaki invariants
is
\[
F(X)=-8\kappa.
\]
\end{cor}

\begin{cor}\label{cor62}  
With the notations as in Theorem~\ref{fud},
if $M$ is the variety of the intersection of
two quadratic polynomials in $CP^4$, then
\[
F(X)=-10(\kappa_1+\kappa_2).
\]
\end{cor}

The first four examples are due to Ding and Tian~\cite{DT}:

\begin{ex}
Let $X_f\subset CP^3$ be the zero locus of a cubic polynomial $f$. Put
$f=z_0z_1^2+z_2z_3(z_2-z_3)$, where $z_0,z_1,z_2,z_3$ are homogeneous
coordinates of $CP^3$. The $X_f$ has a unique quotient singularity at
$p_0=[1,0,0,0]$. This singularity is of the form $C^2/\Gamma$, where
$\Gamma$ is the dihedral subgroup in $SU(2)$ of type $D_4$. One can check
that $X_f$ is a K\"ahler orbifold with $c_1(X)>0$. Let $X$ be the
holomorphic vector field whose real part generates the one parameter
subgroup
$\{diag (1, e^{3t},e^{2t},e^{2t})\}_{t\in R}$ in $SL(4,C)$. Then $X$
restricts to a holomorphic vector field on $X_f$ and has five zeros
$[1,0,0,0]$, $[0,1,0,0]$, $[0,0,1,0]$,$[0,0,0,1]$ and $[0,0,1,1]$. 
\end{ex}

We are going to use three methods to compute the Futaki invariants. 
The first method is the original method in ~\cite{DT}.

We can rewrite the function $f$ near $[1,0,0,0]$ 
in the standard form
\[
f=z_1^2-z_3(z_2^2-4z_3^2).
\]
In~\cite{DPT} we see that there
are standard coverings $C^2\rightarrow C^2/\Gamma$ by
\[
\left\{
\begin{array}{l}
z_1=uv(u^4-v^4)\\
z_2=u^4+v^4\\
z_3=u^2v^2
\end{array}
\right..
\]

If we assume that on the $u-v$ plane, $X=au\frac{\pa}{\pa
u}+bv\frac{\pa}{\pa v}$, then we would have $a=b=\frac 12$.
Since the order of the group $D_4$ is 8, a computation 
using the formula in~\cite[pp 324]{DT} shows
\[
F(X)=\frac 18\cdot \frac{1^3}{\frac 14}
+\frac{(-2)^3}{1}
+3\frac{(-1)^3}{-2}
=-6.
\]

Our second method is a trick which can be generalized to give another
prove of the main theorem of this paper
in the case of hypersurfaces. Suppose that in
$X=au\frac{\pa}{\pa
u}+bv\frac{\pa}{\pa v}$, we don't know what $a$ and $b$ is. By using the
Bott residue formula, we see that
\[
\frac 18\cdot\frac{(a+b)}{ab}
+\frac{(-2)}{1}
+3\frac{(-1)}{-2}
=0
\]
and
\[
\frac 18\cdot\frac{(a+b)^2}{ab}
+\frac{(-2)^2}{1}  
+3\frac{(-1)^2}{-2}  
=3.
\]

Thus we solved $a=b=\frac 12$. The Futaki invariant is obtained.

The last method is to use Corollary~\ref{cor61}, which gives $F(x)=-6$.

\begin{ex}
Let  $f=z_0z_1^2+z_1z_2^2+z_3^3$ and 
$X=\{diag(1,e^{6t},e^{3t},e^{4t})\}$. Then
Corollary~\ref{cor61} gives $F(X)=-18$.
\end{ex}

\begin{ex} 
Let $f=z_0(z_1^2+z_2^2)+z_3^2z_1$ and $X=diag\{(1,e^{2t},e^{2t},e^t)\}$.
Then Corollary~\ref{cor61} gives $F(X)=-2$.
\end{ex}

\begin{ex}
Let $f=z_0(z_1^2+z_2^2)+z_3^3$. and $X=diag\{(1,e^{3t},e^{3t},e^{2t})\}$.
Then Corollary~\ref{cor61} gives
$F(X)=0$.
\end{ex}

The following  examples are given by Jeffres~\cite{JE}: 
let $[z_0,z_1,w,x,z]$ be the general point in $CP^4$.

\begin{ex} Let
\[
\left\{
\begin{array}{l}
f=z_0z_1+w^2+x^2\\
g=z_1l(w,x)+z^2
\end{array},
\right.
\]
where $l(w,x)$ is a linear function of $w,x$. 
Let $X=diag\{1,e^{2t},e^t,e^t,e^{3/2t}\}$
Then $\kappa_1=-1/5$ and $\kappa_2=4/5$. Using Corollary~\ref{cor62}, 
$F(X)=-6$.
\end{ex}

\begin{ex}
\[
\left\{
\begin{array}{l}
f=z_0z_1+z^2\\
g=z_1^2+wx
\end{array}
\right.
\]
and $X=diag\{(1,e^{2t},e^{2t},e^{2t},e^t)\}$
Then $\kappa_1=-4/5$ and $\kappa_2=6/5$. So the Futaki invariant
$F(X)=-4$ using Corollary~\ref{cor62}.
\end{ex}

The following examples are given by Wu~\cite{Wu}.

\begin{ex}
Let $M_0\subset CP^4$ be the zero locus defined by
\[
f=z_0z_1^2+z_1z_2^2+z_3^3+z_4^3=0
\]
and $X=diag (1, e^{6t}, e^{3t}, e^{4t}, e^{4t})$. Then
$\kappa=9/5$ and by Corollary~\ref{fcor}, $F(X)=-36$.
\end{ex}

\begin{ex}
Let $M_0\subset CP^3$ be the zero locus defined by
\[
f=z_3(z_1^2-z_0z_2)+z_2^3=0
\]
and $X=diag (1, e^t, e^{2t}, e^{4t})$. Then
$\kappa=3/4$ and by Corollary~\ref{fcor}, $F(X)=-2$.
\end{ex}

\bibliographystyle{acm}
\bibliography{bib}

\end{document}